\documentstyle[9pt]{amsart}

\title[Eigenvalue gap theorems]{Eigenvalue gap theorems for a class of nonsymmetric elliptic operators on convex domains}
\author{Jon Wolfson}
\address{Department of Mathematics \\ 
                 Michigan State University \\
                 East Lansing, MI 48824}

\date{\today}

\newtheorem{thm}{Theorem}[section]
\newtheorem{lem}[thm]{Lemma}

\newtheorem{prop}[thm]{Proposition}

\theoremstyle{definition}

\newtheorem{defn}{Definition}[section]

\numberwithin{equation}{section}

\renewcommand{\a}{\alpha}
\renewcommand{\b}{\beta}
\renewcommand{\d}{\delta}
\newcommand{\D}{\Delta}
\newcommand{\e}{\varepsilon}
\newcommand{\g}{\gamma}

\renewcommand{\l}{\lambda}
\renewcommand{\L}{\Lambda}
\newcommand{\n}{\nabla}
\newcommand{\p}{\partial}

\newcommand{\s}{\sigma}

\newcommand{\th}{\theta}
\renewcommand{\t}{\tau}
\renewcommand{\O}{\Omega}
\renewcommand{\o}{\omega}
\newcommand{\r}{\rho}
\newcommand{\k}{\kappa}

\def\Pb{\ifmmode{\Bbb P}\else{$\Bbb P$}\fi}
\def\Z{\ifmmode{\Bbb Z}\else{$\Bbb Z$}\fi}
\def\Q{\ifmmode{\Bbb Q}\else{$\Bbb Q$}\fi}
\def\C{\ifmmode{\Bbb C}\else{$\Bbb C$}\fi}
\def\R{\ifmmode{\Bbb R}\else{$\Bbb R$}\fi}
\def\S{\ifmmode{S^2}\else{$S^2$}\fi}

\makeatletter

\newcommand{\Rmnum}[1]{\expandafter\@slowromancap\romannumeral #1@}
\makeatother

\def\dist{\operatorname{dist}}

\def\dist{\operatorname{dist}}

\def\diam{\operatorname{diam}}

\def\Real{\operatorname{Re}}
\def\Im{\operatorname{Im}}

\def\S{\cal S}

\begin{document}

\maketitle

\setcounter{secnumdepth}{1}

\setcounter{section}{-1}

\section{\bf Introduction}

In the remarkable paper [AC] Andrews and Clutterbuck solve the ``gap conjecture'', that is, they show that the difference between the first and second eigenvalues of the laplacian with convex potential on a convex domain in euclidean space is at least $ \frac{3 \pi^2}{D^2}$. Here $D$ is the diameter of the domain. Somewhat later, Lei Ni [N1] reformulated and expanded some of the techniques introduced in [AC]. Taken together, these papers suggest a general  approach to estimating the eigenvalue gap of a large class of linear second-order elliptic operators on convex domains. In this paper we illustrate how this approach may work by estimating the eigenvalue gap of a class of nonsymmetric linear elliptic operators. 

Let $\O$ be a strictly convex open domain in euclidean space with smooth boundary. The operators $L$ we consider have the form:
\begin{equation}
\label{equ:def-of-L-intro}
Lu = \D u - B \cdot \n u - c u
\end{equation}
where $u$ is a $C^2$ function on $\bar{\O}$, $\D$ is the euclidean laplacian, $B$ is a $C^3$ vector valued function on $\bar{\O}$ and $c$ is a $C^2$ scalar function on $\bar{\O}$. 
Let $\b$ denote the one-form dual to $B$. We will require that, for all $x \in \O$, there is a constant $K$ such that:
\begin{equation}
\label{equ:estimate-b-intro}
|d \b(x)| \leq K \dist(x, \p \O).
\end{equation}
Note that this condition is satisfied if, for example, $B$ has compact support in $\O$ or if $B = \n \phi$ for a $C^4$ function $\phi$ on $\O$. We consider the eigenvalue problem:
$$
Lu = -\l u,
$$
where $u$ satisfies Dirichlet boundary conditions. Since $L$ is not symmetric the eigenvalues and eigenfunctions  need not be real. However the principal eigenvalue $\l_0$ is real  with eigenfunction $u_0$ that is positive on $\O$, vanishes on the boundary and satisfies $|\n u_0| \neq 0$ on the boundary. It is also known that for any other eigenvalue $\l$, Re$(\l) > \l_0$ [N2]. We will show that the gap Re$(\l) - \l_0$ can be bounded below by a positive constant $\a$ that depends on the coefficients $B$ and $c$ and on $u_0$. The constant $\a$ is the eigenvalue gap of an associated regular Sturm-Liouville problem on the interval $[-\frac{D}{2}, \frac{D}{2}]$:
$$
w''(s) + \s |s| w(s) = -\mu w(s)
$$
with $w(-\frac{D}{2}) = w(\frac{D}{2})=0$ where $\s$ is a constant depending on $B$ and $c$ and on $u_0$. We note that the spectral gap of a regular Sturm-Liouville problem on an interval is relatively easy to estimate.

In [AC] the eigenvalue gap is bounded below by the eigenvalue gap of the regular Sturm-Liouville problem on the interval $[-\frac{D}{2}, \frac{D}{2}]$:
$$
w'' + \l w = 0.
$$
with $w(-\frac{D}{2}) = w(\frac{D}{2})=0$. Much of the work in [AC] is the determination that this is the {\it associated  Sturm-Liouville problem}, where {\it associated  Sturm-Liouville problem} has a precise technical meaning described in the next section. In this paper the associated Sturm-Liouville problem  is not obviously related to the operator $L$ and is somewhat arbitrary. However the associated Sturm-Liouville problem must satisfy a number of conditions and it is these that are used to determine the particular problem.

\bigskip

The main theorem is:

\begin{thm}
\label{thm:main-intro}
Let $\O$ be a bounded connected strictly convex open domain in $\R^n$ with smooth  boundary $\p \O$ of diameter $D$. Suppose that the operator (\ref{equ:def-of-L-intro}) has  coefficients $b^i$ and $c$ that satisfy $c \geq 0$ and (\ref{equ:estimate-b-intro}). Consider the eigenvalue problem for $L$ on $\O$:
$$
Lu = -\l u,
$$
where $u$ satisfies Dirichlet boundary conditions. Let $\l_0$ be the principal eigenvalue with eigenfunction $u_0$ and $\l$ be any other eigenvalue.
Then there is a constant $\s > 0$ depending on $|| b^j||_{C^2(\bar{\O})}$,  $||\D b^j||_{C^1(\bar{\O})}$, $||c||_{C^2(\bar{\O})}$, $K$ and on  $u_0$ such that:
$$
\mbox{Re} \; (\l) -  \l_0 > \frac{1}{4}  (\mu_1 - \mu_0)  > 0
$$
where $\mu_0 < \mu_1$ are the first two eigenvalues of the Sturm-Liouville problem on $[-\frac{D}{2}, \frac{D}{2}]$:
$$
w''(s) + \s |s| w(s) = -\mu w(s)
$$
with $w(-\frac{D}{2}) = w(\frac{D}{2})=0$.
\end{thm}
The precise nature of the dependence of $\s$ on the geometry of $u_0$ will be described in a later section.

An immediate consequence of our results are the following two applications to the Bakry-Emery Laplacian.
Recall that for a $C^2$ function $\phi$ the Bakry-Emery Laplacian is given by $$\D_\phi = \D - \n \phi \cdot \n.$$ Consider the operator $L$:
$$
Lu = \D_\phi u - c u.
$$
The  Bakry-Emery Laplacian is formally symmetric with respect to the weighed volume form $e^{-\phi} dv$. Therefore the eigenvalue problem:
$$
Lu = -\l u,
$$
with Dirichlet boundary conditions is has real eigenvalues $\l_0 < \l_1 \leq \l_2 \leq \dots$ and real eigenfunctions [E].  
The operator $L$ is of the form (\ref{equ:def-of-L-intro}) with $B = \n \phi$. Therefore $\b = d \phi$ and thus (\ref{equ:estimate-b-intro}) is trivially satisfied. Theorem \ref{thm:main-intro} remains true but the constant $\s$ now depends only on $||\phi||_{C^4(\bar{\O})}$ and $||c||_{C^2(\bar{\O})}$ and, in particular, is independent of the geometry of $u_0$.

We say a function $c$ is {\it $\phi$-convex} if the function $c - \frac{1}{2} \D \phi + \frac{1}{4} |\n \phi|^2 = c  -  \frac{1}{2} \D_{\frac{1}{2} \phi} \phi$ is convex in the usual sense. Note that if $c$ is convex, in the usual sense, and $\D_{\frac{1}{2} \phi} \phi$ is concave, in the usual sense, then $c$ is $\phi$-convex. We show:

\begin{thm}
If $\phi$ is any $C^4$ function and  $c$ is $\phi$-convex  then the spectral gap for the operator $\D_\phi  - c$ on a convex domain $\O$ satisfies:
$$
 { \l}_1 -  { \l}_0 \geq \frac{3 \pi^2}{D^2}.
$$
\end{thm}
In the case that $\phi$ is a constant this is the result of [AC]. 

\bigskip

We are indebted to Lei Ni for introducing us both to his work [N1] and to [AC] and for interesting discussions. In particular, he pointed out that the gap problem for nonsymmetric elliptic operators.  

\bigskip

\section{\bf The Method of Andrews-Clutterbuck}

Suppose that $\O$ is a strictly convex domain in $\R^n$ and $X$ is a vector field on $\O$. A function $\o(s): \R_+ \to \R$ is called a {\it modulus of expansion} for $X$ if for $x,y \in \O$, $x \neq y$
$$
( X(y) - X(x) ) \cdot \frac{ y-x}{|y-x|} \geq 2 \o ( \frac{|y-x|}{2}).
$$
A function $\eta: \R_+ \to \R$ is called a {\it modulus of continuity} for a complex or real valued  function $f$ on $\O$ if for all  $x,y \in \O$
$$
|f(y) - f(x)| \leq 2 \eta(\frac{|y-x|}{2}).
$$

An important result of Andrews-Clutterbuck unifyng these two concepts in proved in [AC]:

\begin{thm}
\label{thm:basicAC}
Let $\O \subset \R^n$ be a strictly convex domain with smooth boundary and with diameter $D$. Let $Y(x,t)$ be a real valued vector field. Let $z(x,t)$ be a smooth, possibly complex valued, solution of 
\begin{equation}
\label{equ:pde-with-drift}
\frac{\p}{\p t} z = \D z - 2Y \cdot \n z.
\end{equation}
with Neumann boundary condition. Suppose that:
\begin{enumerate}
\item $Y( \cdot, t)$ has modulus of expansion $\o(\cdot, t)$ for each $t > 0$, where $\o(s, t): [0, \frac{D}{2}] \times \R_+ \to \R$ is smooth. \\
\item $z(\cdot, 0)$ has modulus of  continuity $\varphi_0$, where $\varphi_0(s): [0, \frac{D}{2}]  \to \R$ is smooth with $\varphi_0(0) =0$ and $\varphi_0'(s) > 0$ on $[0, \frac{D}{2}]$.\\
\item $\varphi(s,t):  [0, \frac{D}{2}] \times \R_+ \to \R$ satisfies:
\begin{enumerate}
\item $\varphi(s,0) = \varphi_0(s)$ on $[0, \frac{D}{2}]$ \\
\item  $\frac{\p \varphi}{\p t} \geq \varphi'' - 2 \o \varphi'$ on   $[0, \frac{D}{2}] \times \R_+$ \\
\item $\varphi'(s,t) > 0$ on $[0, \frac{D}{2}] \times \R_+$ \\
\item $\varphi(0,t) \geq 0$ for each $t \geq 0$.
\end{enumerate}
\end{enumerate}

Then if  $z(x,t)$ is real valued, $\varphi(s,t)$ is a modulus of continuity of $z(x,t)$.  If  $z(x,t)$ is complex valued, $2 \varphi(s,t)$ is a modulus of continuity of $z(x,t)$.
\end{thm}

\begin{pf}
The real case of this theorem is taken directly from [AC]. The case in which $z(x,t)$  is complex valued is needed in this paper and follows from the real case as follows: Since $Y(x,t)$ is real valued, both $\Real(z)$ and $\Im(z)$ satisfy (\ref{equ:pde-with-drift}). Applying the real case to both $\Real(z)$ and $\Im(z)$ implies that each has modulus of continuity $\varphi(s,t)$. The complex case follows.
\end{pf}

We next outline the method of Andrews-Clutterbuck. Let $\O \subset \R^n$ be a strictly convex domain with smooth boundary and with diameter $D$. Let $L$ be the linear, elliptic, symmetric positive operator on scalar functions:
\begin{equation}
\label{equ:elliptic-AC}
L(u) = \D u - c u,
\end{equation}
on $\O$ where $\D u = \sum_{i=1}^n u_{x_i  x_i}$ and  $c \in C^\infty(\O)$  is a  non-negative function. Let  $0 < \l_0 < \l_1 \leq \dots$ be the eigenvalues and  $u_0, u_1, \dots$ the corresponding eigenfunctions for Dirichlet boundary conditions.  Using that $u_0 > 0$ on $\O$ it can be shown that there is a vector field $Y$ so that the function 
$$
z(x,t) = \frac{e^{-\l_1 t} u_1}{e^{-\l_0 t} u_0}
$$
satisfies the heat equation with drift (\ref{equ:pde-with-drift}) with Neumann boundary condition. The vector field $Y$ is the drift velocity. Suppose that $Y$ has modulus of expansion $\o(\cdot, t)$ for each $t > 0$, where $\o(s, t): [0, \frac{D}{2}] \times \R_+ \to \R$ is smooth. 
 
On the interval $[-\frac{D}{2}, \frac{D}{2}]$ consider the linear  second order differential operator $\tilde{L}(w) = w'' -  \tilde{c} w$ where $\tilde{c}$ is a smooth function on $[-\frac{D}{2}, \frac{D}{2}]$. Consider the Sturm-Liouville problem:
\begin{eqnarray}
\label{equ:Sturm-Liouville}
&&\tilde{L}(w)  =  - \mu w\\ 
&&w(-\frac{D}{2})  = w(\frac{D}{2}) = 0 \nonumber
\end{eqnarray}
with eigenvalues $0 < \mu_0 < \mu_1 \leq \dots$ and corresponding eigenfunctions $w_0, w_1, \dots$. Set 
$$
\varphi (s, t)= \frac{e^{-\mu_1 t} w_1}{e^{-\mu_0 t} w_0}
$$
for $s \in [0, \frac{D}{2}]$ and $t \geq 0$. Then
\begin{equation}
\frac{\p \varphi}{\p t} = \varphi'' - 2 \tilde{\o} \varphi' \;\; \mbox{on} \;\;[0, \frac{D}{2}] \times \R_+
\end{equation}
Here the function $\tilde{\o}(s,t)$ is also called  the drift velocity. It is not difficult to verify by direct computation that:
\begin{equation}
 \tilde{\o} = -(\log w_0)',
\end{equation}
where $w_0$ is the first eigenfunction of (\ref{equ:Sturm-Liouville}). Suppose that the operators $L$ and $\tilde{L}$ satisfy the condition that the potential function $\tilde{c}$ is even and that the potential function $c$ is more convex then  $\tilde{c}$ in the sense that for any $x \neq y$ in $\O$:
\begin{equation}
(\n c (y) - \n c(x)) \cdot \frac{(y-x)}{|y-x|} \geq 2 \tilde{c}'( \frac{|y-x|}{2}).
\end{equation}
Under this assumption [AC] prove that:
\begin{enumerate}
\item $\tilde{\o} = \o$ \\
\item $\varphi'(s,t) > 0$ on $[0, \frac{D}{2}] \times \R_+$ \\
\item $\varphi(0,t) \geq 0$ for each $t \geq 0$.
\end{enumerate}
By Theorem \ref{thm:basicAC}, $\varphi(s,t)$ is a modulus of continuity of $z(x,t)$. Hence, for each $t \geq 0$
$$
\mbox{osc} \; z(\cdot, t) \leq 2 \; \mbox{sup} \; \{ \varphi(s,t) : s \in [0, \frac{D}{2}] \}.
$$
From this they derive that there is a constant $C$ such that for each $t \geq 0$
$$
e^{-(\l_1-\l_0) t} \;  \mbox{osc}\; \frac{u_1}{u_0} \leq 2C e^{-(\mu_1-\mu_0) t}.
$$
Hence,
\begin{equation}
\label{equ:gap-inequality}
\l_1-\l_0 \geq \mu_1-\mu_0.
\end{equation}

\bigskip

In this paper we will employ a variation of this argument.  As above let  $\O$ be a bounded connected strictly convex open domain in $\R^n$ with smooth  boundary $\p \O$. Consider the uniformly elliptic operator:
$$
L u = \sum_{i=1}^n  u_{x_i x_i} - \sum_{i=1}^n b^i u_{x_i} - cu.
$$
with
$$
|d \b(x)| \leq K \dist(x, \p \O),
$$
where $\b$ is the one form dual to $(b^j)$. For Dirichlet boundary conditions let $\l_0$ be the principal eigenvalue with corresponding eigenfunction $u_0$ and let $\l \in \C$ be any other eigenvalue with eigenfunction $u$.  We will show that there is a vector field $Y$ so that the function 
$$
z(x,t) = \frac{e^{-\l t} u}{e^{-\l_0 t} u_0}
$$
satisfies the heat equation with drift (\ref{equ:pde-with-drift}) with Neumann boundary condition. Moreover $Y$ has modulus of expansion $\o(\cdot, t)$ for each $t > 0$. 

We next consider the Sturm-Liouville problem on $[-\frac{D}{2}, \frac{D}{2}]$:
\begin{equation}
\label{equ:SLproblem}
w''(s) + \s |s| w(s) = -\mu w(s)
\end{equation}
with $w(-\frac{D}{2}) = w(\frac{D}{2})=0$.  Let $\mu_0 < \mu_1$ be the first two eigenvalues with corresponding eigenfunctions $w_0$ and $w_1$. We show that for suitable choice of $\s$ and a scaling factor $\eta > \frac{1}{2}$ the functions $e^{-\mu_0 \eta^2 t} w_0$ and $e^{-\mu_1 \eta^2 t} w_1$ can be used to define a function $\varphi$ on $[0, \frac{D}{2}] \times \R_+$ such that:
$$
\frac{\p \varphi}{\p t} \geq \varphi'' - 2 \o \varphi'
$$ 
on   $[0, \frac{D}{2}] \times \R_+$, where $\o$ is the modulus of expansion of $Y$. The other conditions needed for Theorem  \ref{thm:basicAC} can also be verified for this choice of $\varphi$. The
argument proceeds as in [AC] to conclude that:
\begin{equation}
\label{equ:gap-inequality2}
\mbox{Re} \; (\l) -  \l_0 \geq \eta^2 (\mu_1 - \mu_0) >   \tfrac{1}{4} (\mu_1 - \mu_0) > 0
\end{equation}
 We will call the Sturm-Liouville problem (\ref{equ:SLproblem}) an {\it associated Sturm-Liouville problem} for $L$.
The difference $ \mu_1-\mu_0$ can be computed (or estimated) by applying ode techniques to the associated Sturm-Liouville problem, thus providing a lower bound on the eigenvalue gap of the operator $L$. In the case that the one-dimensional limit of the operator $L$ coincides with an associated Sturm-Liouville problem the lower bound (\ref{equ:gap-inequality2}) is sharp.

Using our nomenclature the regular Sturm-Liouville problem on the interval $[-\frac{D}{2}, \frac{D}{2}]$:
$$
w'' + \l w = 0.
$$
with $w(-\frac{D}{2}) = w(\frac{D}{2})=0$ is an associated Sturm-Liouville problem for (\ref{equ:elliptic-AC}) under the condition that $c$ is convex. This is the technical meaning of ``associated'' referred to in the introduction.

\bigskip

\section{\bf The gap theorem for nonsymmetric elliptic operators}

Let $\O$ be a bounded connected strictly convex open domain in $\R^n$ with smooth  boundary $\p \O$. Assume that $ b^i \in C^{3}(\bar \O)$ and that $c \in C^{2}(\bar \O)$ and consider the uniformly elliptic operator:
\begin{equation}
\label{equ:def-of-L}
L u = \sum_{i=1}^n  u_{x_i x_i} - \sum_{i=1}^n b^i u_{x_i} - cu.
\end{equation}
Consider the one-form $\beta =  \sum_{i=1}^n b^i dx_i$ dual to the vector field $ \sum_{i=1}^n b^i \frac{\p}{\p x_i}$ and suppose that there is a constant $K > 0$ such that:
\begin{equation}
\label{equ:estimate-b}
| d \b (x)| \leq K \dist(x , \p \O),
\end{equation}
for all $x \in \O$. We will, throughout this paper, assume that the coefficients of $L$ satisfy (\ref{equ:estimate-b}). Note that we can add a constant to the operator $L$ without changing the spectral gap of $L$. Therefore we can, without loss of generality, assume that $c > 0$ on $\O$. The operator $L$ is not symmetric and the eigenvalues for Dirichlet boundary conditions may not all be real. However the following theorem quoted from [E] shows that, in part, the situation resembles the symmetric case.

\begin{thm}

\begin{enumerate}
\item The principle eigenvalue $\l_0$ of $L$ on $\O$ with zero boundary conditions is real and simple.\\

\item If $\l \in \C$ is any other eigenvalue of $L$ then
$$
\mbox{Re} \; \l \geq \l_0
$$\
\item The eigenfunction $u_0$ corresponding to $\l_0$ is positive in $\O$.\\
\end{enumerate}
\end{thm}

Recently Lei Ni [N2] has improved this result showing that if $\l \in \C$ is any other eigenvalue of $L$  then the strict inequality $\mbox{Re} \; \l > \l_0$ holds.

\bigskip

We will prove the gap theorem Theorem \ref{thm:main-intro} for the eigenvalues of the operator (\ref{equ:def-of-L}). We begin with the following proposition. 

\begin{prop}
\label{prop:bound-on-beta}
Let $\O$ be a bounded connected strictly convex open domain in $\R^n$ with smooth  boundary $\p \O$. Suppose that $\l_0$ is the principal eigenvalue of the operator (\ref{equ:def-of-L}) with eigenfunction 
$u_0$.  Then, there is a constant $\k$ such that:
\begin{equation}
\sup_{x \in \O} \frac{|\n u_0(x)|}{u_0(x)} | d \b (x)|  \leq \k.
\end{equation}
\end{prop}

\begin{pf}
The principal  eigenfunction $u_0$ satisfies $u_0 > 0$ on $\O$, $u_0{_{|_{\p \O}}} = 0$, and $\frac{\p u_0}{\p \nu}_{|_{\p \O}} < 0$. There is a constant $A > 0$ such that: 
\begin{equation}
\label{equ:estimate-A}
| u_0|_{C^1(\bar{\O})} \leq A
\end{equation}
Set $\O_\d = \{ x \in \O : u_0(x) \geq \d \}$. Since $|\n u_0| > 0$ on $\p \O$, there exist constants $\th_0 > 0$ and $\d_0> 0$ such that on $\O \setminus \O_{\d_0}$:
\begin{equation}
\label{equ:estimate-u}
|\n u_0| \geq \th_0.
\end{equation}
It follows that all the critical points of $u_0$ occur in $\O_{\d_0}$. Set $a = \min_{\O_{\d_0}} u_0$. Then for $x \in \O_{\d_0}$ we have: 
$$
 \frac{|\n u_0(x)|}{u_0(x)} | d \b (x)|  \leq \frac{A}{a}K D, 
$$
where $D = \diam(\O)$.
Choosing ${\d_0}$ smaller, if necessary, we can suppose that the line $\g$ joining $x \in \O \setminus  \O_{\d_0}$ to its nearest boundary point $y$ lies entirely in $\O \setminus  \O_{\d_0}$. Parameterizing 
$\g$ by its arc length with $\g(0) = y$ and $\g(\ell) = x$ we can also suppose that $-\n u_0(\g(s)) \cdot \g'(s) \geq \frac{\th_0}{2}$ along $\g$. Thus,
$$
u_0(x) = u_0(x) - u_0(y) = \int_0^\ell -\n u_0(\g(s)) \cdot \g'(s) ds \geq  \frac{ \th_0 \ell}{2},
$$
where $\ell$ is the length of $\g$. Hence for $x \in \O \setminus  \O_{\d_0}$ we have:
$$
 \frac{|\n u_0(x)|}{u_0(x)} | d \b (x)|  \leq \frac{2AK}{\th_0 \ell}   \dist(x, \p \O) \leq \frac{2A K}{ \th_0}.
$$
The result follows.
\end{pf}

The dependence of the constant $\k$ on the geometry of $u_0$ can be described as follows: Denote the set of critical points of $u_0$ in the interior of $\O$ by $S$. Set $\d_0 = \frac{1}{2} \inf_{x \in S}  u_0(x)$. Using the notation of the proof, we have that $\inf_{x \in \O \setminus \O_{\d_0}} |\n u_0 (x)| =  \th_0 > 0$.  Then $\k$ depends on $\d_0, \th_0, |\n u_0|_{C^1(\bar{\O})}$  and $K$.

\bigskip

\subsection{The drift velocity of $L$}

Let $u_0$ be the principal eigenfunction of $L$ with eigenvalue $\l_0$ and set $u_0(x,t) = e^{-\l_0 t} u_0(x)$.
Let $u$ be any other eigenfunction with eigenvalue $\l$ and set $u_1(x,t) = e^{-\l t} u(x)$. The following proposition is adapted from [AC].

\begin{prop}
Let $\O$ be a bounded strictly convex domain  with smooth boundary in $\R^n$.
Let $u_0$ and $u_1$ be two smooth solutions of the parabolic equation:
\begin{eqnarray*}
\frac{\p u}{\p t} & =& L( u ) \;\; \mbox{on } \;\; \O \times \R_+\\
u &=& 0 \;\; \mbox{on } \;\; \p \O \times \R_+
\end{eqnarray*}
with $u_0$ is positive on the interior of $\O$.
Let $z(x, t) = \frac{u_1(x,t)}{u_0(x, t)}$ and  $$Y^j (x) = -  \n_{x_j} (\log u_0)(x) + \frac{b^j}{2}(x).$$ Then $z$ is smooth  on $\O \times \R_+$  and  satisfies the Neumann heat equation with drift:
\begin{eqnarray}
\label{equ:drift1}
\frac{\p z}{\p t} & = & \D z - 2 Y \cdot \n z \;\; \mbox{on } \;\; \O \times \R_+\\
\n_\nu z &=& 0 \;\;   \mbox{on } \;\; \p \O \times \R_+
\end{eqnarray}
\end{prop}

\begin{pf}
The proof is essentially given in [AC] Proposition 3.1 (or [Y] Lemma 1.1, [SWYY] Appendix A). Both $u_0$ and $u_1$ are smooth on $\bar{\O} \times [0, \infty)$ and $u_0$ has negative derivative in the direction of the inward pointing unit normal. By the argument of [SWYY] $z$ extends to $\bar{\O}$ as a smooth function and therefore $\frac{\p z}{\p t}, \D z$ and $\n z$ are smooth and bounded on $\bar{\O}$.
By direct computation:
\begin{eqnarray} \nonumber
\frac{\p z}{\p t} &=&\frac{\p }{\p t} (\frac{u_1}{u_0})\\ 
\label{eqn:drift}
&=& \D z +  ( 2 \n \log u_0 - B) \cdot \n z
\end{eqnarray}
On $\p \O$, $\n u_0 = -k \nu$ with $k > 0$ and $u_0 =0$. Therefore by (\ref{eqn:drift}), $\n_\nu z = 0$.
\end{pf}

\noindent $Y$ is called  the {\it drift velocity} of $z$.

\begin{lem}
\label{lem:laplacian-on-Y}
Set $B = (b^j)$, $U^{ij}= \n_{x_j} b^i - \n_{x_i} b^j$ and $$V^j = V^j(c, B) =  \big( \n_{x_j} c + \tfrac{1}{4}  \n_{x_j} (|B|^2) - \tfrac{1}{2} \D b^j \big).$$  Then
\begin{equation}
\label{equ:laplacian-of-Y}
\D Y = 2 \n_Y Y - Y \cdot U - V
\end{equation}
\end{lem}

\begin{pf}
To begin we compute $\n_Y Y$. 
\begin{eqnarray} \nonumber
\n_Y Y & = & \sum_i (-\n_{x_i} \log u_0 + \frac{b^i}{2}) \n_{x_i} Y^j \\ \nonumber
&=&  \sum_i (-\n_{x_i} \log u_0 + \frac{b^i}{2})  ( -\n_{x_j} \n_{x_i} \log u_0 + \frac{1}{2} \n_{x_i} b^j )\\ \nonumber
\label{equ:gradient-of-Y}
&=& \sum_i \big( (\n_{x_i} \log u_0)( \n_{x_j} \n_{x_i} \log u_0) -   \frac{b^i}{2} (\n_{x_j} \n_{x_i} \log u_0) \\ 
&-&  \frac{1}{2}   (\n_{x_i} \log u_0) \n_{x_i} b^j + \frac{1}{4} b^i \n_{x_i} b^j \big)
\end{eqnarray}
Also we will need:
\begin{eqnarray*}
\sum_i  \n_{x_i}  \n_{x_i} \log u_0 &=& \sum_i  \big( \frac{1}{u_0}( \n_{x_i}  \n_{x_i} u_0) - (\n_{x_i} \log u_0)  (\n_{x_i} \log u_0) \big)\\
&=& \sum_i \big( b^i \n_{x_i} \log u_0 + c - \l_0 - (\n_{x_i} \log u_0) ( \n_{x_i} \log u_0) \big)
\end{eqnarray*}
Computing $\D (- \n_{x_j} \log u_0) = - \n_{x_j} ( \D \log u_0)$ we have, 
\begin{eqnarray} \nonumber
\label{equ:laplacian-of-log}
\sum_i -\n_{x_j} (\n_{x_i}  \n_{x_i} \log u_0) &=&\sum_i \big( - \n_{x_j} b^i ( \n_{x_i} \log u_0) - b^i (\n_{x_j}\n_{x_i} \log u_0) \\
&-& \n_{x_j} c + 2 ( \n_{x_j} \n_{x_i} \log u_0) ( \n_{x_i} \log u_0) \big)
\end{eqnarray}
Combining (\ref{equ:gradient-of-Y}) and (\ref{equ:laplacian-of-log}) we have:
\begin{eqnarray*}
\sum_i -\n_{x_j} (\n_{x_i}  \n_{x_i} \log u_0) &=& 2 \n_Y Y + \sum_i  ( \n_{x_i} b^j -  \n_{x_j} b^i)( \n_{x_i} \log u_0)  \\
&-& \n_{x_j} c    - \frac{1}{2} \sum_i b^i \n_{x_i} b^j  \\
&=& 2 \n_Y Y + \sum_i  ( \n_{x_i} b^j -  \n_{x_j} b^i)( \n_{x_i} \log u_0 - \frac{b^i}{2})  \\
&-&( \n_{x_j} b^i)  \frac{b^i}{2} - \n_{x_j} c     \\
\end{eqnarray*}
Therefore
\begin{eqnarray*}
\D Y &=& \D (- \n_{x_j} \log u_0 + \frac{b^j}{2} )\\
&=& 2 \n_Y Y - \sum_i  ( \n_{x_i} b^j -  \n_{x_j} b^i) Y^i -  \frac{1}{4} \n_{x_j} (\sum_i b^i b^i) - \n_{x_j} c + \D(\frac{b^j}{2}). \\
\end{eqnarray*}
\end{pf}

\begin{lem}
\label{lem:estimate-on-Y}
On $\O$ there is a constant $\L$ depending on $\k$ and on $\sup_j || b^j ||_{C^1(\bar{\O})}$ such that:
\begin{equation}
\sup_{x \in \O} |Y(x)| | U(x)| \leq \L.
\end{equation}
\end{lem}

\begin{pf}
This follows easily from Proposition \ref{prop:bound-on-beta}.
\end{pf}

We suppose that there is a function $\t: \R_+ \to \R$ such that on $\O$:
\begin{equation}
\label{equ:estimate-on-V}
\big( V(y) - V(x) \big) \cdot \frac{y-x}{|y-x|} \geq 2\t( \frac{|y-x|}{2})
\end{equation}
Let $\psi(s): [0, \frac{D}{2}) \to \R$ be a $C^2$ function which satisfies for each $s \in [0, \frac{D}{2})$:
\begin{equation}
\label{equ:psi-neg}
\psi(0) \geq 0, \;\;\; \psi' (s) < 0
\end{equation}
\begin{equation}
\label{equ:psi-eq1}
2\L - 2\t(s) + 2 \psi''(s) \leq - 4\psi'(s) \psi(s).
\end{equation}

\bigskip

The following theorem is motivated by a similar result in [AC] and [N1]. It is the first step in deriving a modulus of expansion for $Y$.

\begin{thm}
\label{thm:estimate-on-C}
Suppose that $Y$ satisfies (\ref{equ:laplacian-of-Y}) on $\O$. Let $\psi$ be the function defined above. Then 
$$
{\cal C}(x,y) = \big( Y(y) - Y(x) \big) \cdot \frac{y-x}{|y-x|} + 2 \psi \bigg( \frac{|y-x|}{2} \bigg)
$$
can not attain a negative minimum in the interior of $\O$.
\end{thm}

\begin{pf} We argue by contradiction and assume that  at $(x_0, y_0)$, ${\cal C}(x,y)$ attains a negative minimum. Clearly $x_0 \neq y_0$ since  ${\cal C}(x,x) \geq 0$.
Following [AC] and [N1] we choose a local orthonormal frame at $x_0$, denoted $\{ e_1, \dots, e_n \}$, with $e_n = \frac{y_0-x_0}{|y_0 - x_0|}$ and parallel translate this frame along the line interval joining $x_0$ to $y_0$. Then at $(x_0, y_0)$ we derive that:
\begin{eqnarray} \nonumber
0 &=& \frac{\p}{\p s} {\cal C}(x + s e_i,y)_{|_{s=0}} \;\;\;\; \mbox{for} \; 1 \leq i \leq n-1, \\
&=& -\n_{e_i} Y(x) \cdot  \frac{y-x}{|y-x|}  - \frac{Y(y) -Y(x)}{|y-x|} \cdot e_i. 
\end{eqnarray}
\begin{eqnarray} \nonumber
0 &=& \frac{\p}{\p s} {\cal C}(x ,y + s e_i)_{|_{s=0}} \;\;\;\; \mbox{for} \; 1 \leq i \leq n-1, \\
&=& \n_{e_i} Y(y) \cdot  \frac{y-x}{|y-x|}  + \frac{Y(y) -Y(x)}{|y-x|} \cdot e_i. 
\end{eqnarray}
\begin{eqnarray} \nonumber
0 &=& \frac{\p}{\p s} {\cal C}(x + s e_n,y)_{|_{s=0}}, \\
&=& -\n_{e_n} Y(x) \cdot  \frac{y-x}{|y-x|}  - \psi' \bigg( \frac{|y-x|}{2} \bigg).
\end{eqnarray}
\begin{eqnarray} \nonumber
0 &=& \frac{\p}{\p s} {\cal C}(x ,y + s e_n)_{|_{s=0}}, \\
&=& \n_{e_n} Y(y) \cdot  \frac{y-x}{|y-x|}  + \psi' \bigg( \frac{|y-x|}{2} \bigg).
\end{eqnarray}
Let $E_i = e_i \oplus e_i \in T_{(x_0, y_0)} \R^n \times \R^n$ for $1 \leq i \leq n-1$ and $E_n = e_n \oplus (-e_n)$. Since ${\cal C}(x,y)$ attains its minimum at $(x_0, y_0)$ we have:
$$
\n^2_{E_i E_i} {\cal C}_{|_{(x_0, y_0)}} \geq 0, \; \mbox{for} \; 1 \leq i \leq n.
$$
Along the path $(x + s e_i, y + s e_i)$ for $1 \leq i \leq n-1$ we note that $y-x$ is constant. Thus, computing as in [AC], we see that at $(x_0, y_0)$:
\begin{equation}
0 \leq \n^2_{E_i E_i} {\cal C} =  \big( \n^2_{e_i e_i}Y(y) - \n^2_{e_i e_i}Y(x) \big) \cdot \frac{y-x}{|y-x|}  \; \mbox{for} \; 1 \leq i \leq n-1.
\end{equation}
Along the path $(x + s e_n, y - s e_n)$, $ \frac{y-x}{|y-x|}$ is constant, $\frac{d}{ds} |y-x| = -2$ and $\frac{d^2}{ds^2} |y-x| = 0$. Again, computing as in [AC], we derive that at $(x_0, y_0)$:
\begin{equation}
0 \leq \n^2_{E_n E_n} {\cal C}  =  \big( \n^2_{e_n e_n}Y(y) - \n^2_{e_n e_n}Y(x) \big) \cdot \frac{y-x}{|y-x|} + 2 \psi''
\end{equation}
Using Lemma \ref{lem:laplacian-on-Y} we have that at $(x_0, y_0)$:
\begin{eqnarray*}
&&\sum_{i=1}^n \n^2_{E_i E_i} {\cal C} = \big( \D Y(y) - \D Y(x) \big) \cdot \frac{y-x}{|y-x|} + 2 \psi'' \\
&=& 2 \big( \n_{Y(y)} Y(y) - \n_{Y(x)} Y(x) \big) \cdot \frac{y-x}{|y-x|} + \big(- Y(y) \cdot U(y) +  Y(x) \cdot U(x) \big) \cdot \frac{y-x}{|y-x|} \\
&-& \big(V(y) -V(x) \big) \cdot \frac{y-x}{|y-x|} +  2 \psi''  \\
\end{eqnarray*}
As in [AC] and [N1] note that at $(x_0, y_0)$:
\begin{eqnarray*}
 \n_{Y(y)} Y(y)  \cdot \frac{y-x}{|y-x|} &=& \langle \n_{Y(y)} Y(y), e_n \rangle \\
 &=& \sum_{i=1}^n \langle Y(y), e_i \rangle \langle \n_{e_i} Y(y), e_n \rangle \\
 &=& \frac{-1}{|y-x|}\sum_{i=1}^{n-1} \langle Y(y), e_i \rangle \langle  Y(y) - Y(x) , e_i \rangle  - \psi' Y(y) \cdot \frac{y-x}{|y-x|} \\ 
 \end{eqnarray*}
By Lemma \ref{lem:estimate-on-Y}:
\begin{equation}
 \big| \big(- Y(y) \cdot U(y) +  Y(x) \cdot U(x) \big)  \cdot \frac{y-x}{|y-x|} \big| \leq 2 \L
 \end{equation}

Putting these inequalities together we have at $(x_0, y_0)$:
\begin{eqnarray*}
\sum_{j=1}^n \n^2_{E_j E_j} {\cal C} &\leq& \frac{-2}{|y-x|}\sum_{j=1}^{n-1} \langle  Y(y) - Y(x) , e_i \rangle^2\\
 & -  & 2\psi' (s) (Y(y) - Y(x))  \cdot \frac{y-x}{|y-x|} + 2\L - 2\t(s) + 2 \psi'' (s)\\ 
 & \leq & -2\psi' (s) (Y(y) - Y(x))  \cdot \frac{y-x}{|y-x|}  - 4\psi'(s)  \psi(s)\\
 & \leq & -2\psi' (s) {\cal C}(x,y)\\
 & < & 0
\end{eqnarray*}
Here $s = \frac{|y-x|}{2}$. The second inequality uses (\ref{equ:psi-eq1}). 
The conclusion contradicts the assumption that $(x_0, y_0)$ is a minimum point of ${\cal C}(x,y)$.
\end{pf}

\bigskip

\subsection{Boundary asympotics of ${\cal C}(x,y)$}

The next step in the derivation of a modulus of expansion of $Y$ is to
study the boundary behavior of the function ${\cal C}(x,y)$. This analysis is similar to that done in  [N1] and [AC], though we must modify the argument to our situation. In particular, unlike the situation in [N1] and [AC] we do not have available the log convexity of the first eigenfunction $u_0$. We continue under the assumption that there is a $C^2$ function $\psi$ on $[0, \frac{D}{2})$ that satisfies (\ref{equ:psi-neg}) and (\ref{equ:psi-eq1}).  Note first that $\psi(s)$ is not defined at $s = \frac{D}{2}$. To rectify this, for fixed $D' > D$ consider $\psi$ to be a solution of (\ref{equ:psi-eq1}) on $[0, \frac{D'}{2})$. Then $\psi$ is uniformly continuous on $[0, \frac{D}{2}]$. 

Set $X = - \n \log u_0$ so that $Y = X + \frac{1}{2} B$. Then:
$$
{\cal C}(x,y) = \big( X(y) - X(x) \big) \cdot \frac{y-x}{|y-x|} + \frac{1}{2} \big( B(y) - B(x) \big) \cdot \frac{y-x}{|y-x|}+ 2 \psi \bigg( \frac{|y-x|}{2} \bigg)
$$
Given $\e > 0$ we want to show that ${\cal C}(x,y) \geq -\e$ on $\O \times \O$. To begin we note that $\frac{1}{2} \big( B(y) - B(x) \big) \cdot \frac{y-x}{|y-x|}+ 2 \psi ( \frac{|y-x|}{2} ) \geq 0$ on the diagonal  $\D = \{(x,x) : x \in \O \}$. Use the uniform continuity of $\psi$  on $[0, \frac{D}{2}]$ to find a neighborhood of the diagonal $\D_\eta = \{ (x,y) \in \O \times \O: |x-y| < \eta \}$ on which:
\begin{equation}
\label{equ:lower-bound}
\frac{1}{2} \big( B(y) - B(x) \big) \cdot \frac{y-x}{|y-x|}+ 2 \psi \bigg( \frac{|y-x|}{2} \bigg) \geq -\frac{\e}{2}
\end{equation}
The main work of this subsection  involves the term $\big( X(y) - X(x) \big) \cdot \frac{y-x}{|y-x|}$.
Here $u_0$ is the first eigenfunction of the operator $L$ on $\O$ with Dirichlet boundary conditions. It satisfies $u_0 > 0$ on $\O$, $u_0{_{|_{\p \O}}} = 0$, and $\frac{\p u_0}{\p \nu}_{|_{\p \O}} < 0$. There is a constant $A > 0$ such that: 
\begin{equation}
\label{equ:estimate-A}
| u_0|_{C^2(\bar{\O})} \leq \frac{A}{2}
\end{equation}

A neighborhood of the diagonal  must be treated separately. This is because on the diagonal $\big( X(y) - X(x) \big) \cdot \frac{y-x}{|y-x|} = 0$ and therefore its behavior as $x \to \p \O$ differs from its behavior away from the diagonal. We require two results. Set $\O_\d = \{ x \in \O : u_0(x) \geq \d \}$. The first result states that for $\d$ and $\eta$ sufficiently small, if $x, y \in (\O \setminus \O_\d) \cap \D_\eta$ then
$$
\big( X(y) - X(x) \big) \cdot \frac{y-x}{|y-x|} \geq 0.
$$
Thus there is an $\eta_0 \leq \eta$ so that on $\D_{\eta_0}$,  $\big( X(y) - X(x) \big) \cdot \frac{y-x}{|y-x|} \geq -\frac{\e}{2}$. The second result shows that for $\d$ sufficiently small,  if $x, y \in (\O \setminus \O_\d)$ and $|x-y| > \eta_0$ then $\big( X(y) - X(x) \big) \cdot \frac{y-x}{|y-x|}$ is large. In fact, on this set $\big( X(y) - X(x) \big) \cdot \frac{y-x}{|y-x|} \to \infty$ as $\d \to 0$.

To prove both results we study $\big( X(y) - X(x) \big) \cdot \frac{y-x}{|y-x|}$ on $\O \times \O$. 
Since $|\n u_0| > 0$ on $\p \O$, there exist constants $\th_0 > 0$ and $\d' > 0$ such that on $\O \setminus \O_{\d'}$:
\begin{equation}
\label{equ:estimate-u}
|\n u_0| \geq \th_0
\end{equation}
By the implicit function theorem this implies that for each $\d \leq \d'$ the set $\p \O_{\d}$ is a smooth hypersurface. By the convexity of $\O$ it follows that for $\d'$ sufficiently small there is
a constant $\th_1 > 0$ such that for each $\d \leq \d'$  the second fundamental form $II(\cdot, \cdot)$ of the hypersurface $\p \O_\d$ satisfies the inequality:
\begin{equation}
\label{equ:estimate-II}
{\Rmnum{2}}(\cdot, \cdot) \geq \th_1 I(\cdot, \cdot),
\end{equation}
where $I(\cdot, \cdot)$ is the metric on $\p \O_\d$. On $\p \O_\d$ the second fundamental form is given by:
\begin{equation}
\label{equ:define-II}
II(\cdot, \cdot) = \frac{\n^2 u_0(\cdot, \cdot)}{|\n u_0|},
\end{equation}
since $\p \O_\d$ is a level set of $u_0$. For use below we set $C_1 = \frac{1}{2} \big( \frac{ A^2}{\th_0 \th_1} + A \big)$ and $\d'' = \min\{\d', \frac{\th_0^2}{4C_1} \}$. 

\begin{lem}
\label{lem:first-tech}
There is a $\bar{\d} < \frac{1}{2} \d''$ and an $\eta_0 < \eta$ such that for $x,y \in (\O \setminus \O_{\bar{\d}}) \cap \D_{\eta_0}$, $x \neq y$:
$$
\big( X(y) - X(x) \big) \cdot \frac{y-x}{|y-x|} \geq 0
$$
\end{lem}

\begin{pf}
Choose $\bar{\d} < \frac{1}{2} \d''$ and $\eta_0 < \eta$ such that for any two points  $x,y \in (\O \setminus \O_{\bar{\d}}) \cap \D_{\eta_0}$, $x \neq y$, the line segment joining $x$ to $y$ lies in $\O \setminus \O_{\d''}$. Let $\g(s)$ denote this line segment parameterized by arc length. Without loss of generality we can suppose that $u_0(x) \leq u_0(y) \leq \bar{\d}$. Set $\d = u_0(x)$. Then,
\begin{eqnarray} \nonumber
(X(y) - X(x) ) \cdot \frac{y-x}{|y-x|} & = & \langle X(\g(s)), \g'(s) \rangle|_0^{|y-x|} \\ \nonumber
\label{equ:integral-est}
& = & \int_0^{|y-x|} \frac{d}{ds} ( \langle X(\g(s)), \g'(s) \rangle) ds\\ 
& = & \int_0^{|y-x|} \n^2(-\log u_0)( \g'(s), \g'(s) ) ds\\ \nonumber
\end{eqnarray}
Set $\g'(s) = W$. At each point $\g(s)$, decompose $W$ into a component, $W^{\top}$, tangent to $T_{\g(s)} \O_{u_0(\g(s))}$ and a normal component $W^\perp$ with respect to the inward pointing normal $-\nu_{\g(s)}$. We have:
\begin{eqnarray*}
\n^2 u_0(W,W) & = & \n^2 u_0(W^{\top},W^{\top}) + 2\n^2 u_0(W^{\top},W^\perp)+ \n^2 u_0(W^\perp,W^\perp)\\
& \leq & -|\n u_0| II(W^{\top},W^{\top}) + A |W^{\top}||W^\perp| + \frac{A}{2} |W^\perp|^2\\
& \leq & - \th_0 \th_1|W^{\top}|^2 + A |W^{\top}||W^\perp| +  \frac{A}{2} |W^\perp|^2\\
& \leq & - \frac{\th_0 \th_1}{2} |W^{\top}|^2 + C_1 |W^\perp|^2\\
\end{eqnarray*}
In the last inequality we have used:
$$
A |W^{\top}||W^\perp| = (\th_0 \th_1)^{\frac{1}{2}} |W^{\top}| \frac{A}{(\th_0 \th_1)^{\frac{1}{2}}} |W^\perp| \leq \tfrac{1}{2} \big( \th_0 \th_1|W^{\top}|^2 + \frac{A^2}{\th_0 \th_1} |W^\perp|^2 \big)
$$
Since $u_0(x) = \d \leq \bar{\d}$, we  let $k$ be the integer such that
$$
2^k \d \leq \bar{\d} < 2^{k+1} \d < \d''.
$$
For $j= 1, \dots, k+1$ set $\d_j = 2^j \d$. Then for $\d_{j-1} \leq u_0 \leq \d_{j}$ we have:
\begin{eqnarray*}
\n^2 \log u_0(W,W) & = & \frac{\n^2 u_0(W,W)}{u_0}  - \frac{|\n u_0|^2}{u_0^2} |W^\perp|^2 \\
& \leq &- \frac{\th_0 \th_1}{2 u_0}|W^{\top}|^2 + \frac{C_1}{ u_0}  |W^\perp|^2 -  \frac{\th_0^2}{u_0^2}|W^\perp|^2 \\
& \leq &- \frac{\th_0 \th_1}{2 \d_{j}}|W^{\top}|^2 + \frac{2C_1}{\d_{j}}  |W^\perp|^2 -  \frac{\th_0^2}{{\d_{j}}^2}|W^\perp|^2 \\
& \leq &- \frac{\th_0 \th_1}{2 \d_{j}}|W^{\top}|^2  -  \frac{\th_0^2}{2{\d_{j}}^2}|W^\perp|^2 \\
\end{eqnarray*}
where the final inequality uses the definition of $\d''$ and that $\d_j < \d''$. From this inequality  it follows immediately that:
\begin{equation}
\label{equ:integral-positive}
\int_{{\g}} \n^2(-\log u_0)( \g'(s), \g'(s) ) ds \geq 0
\end{equation}
The result follows from (\ref{equ:integral-est}).
\end{pf}

\bigskip

It follows from the lemma that there is an $\eta_1 \leq \eta_0$ so that on $\D_{\eta_1}$,  $\big( X(y) - X(x) \big) \cdot \frac{y-x}{|y-x|} \geq -\frac{\e}{2}$. We next study $\big( X(y) - X(x) \big) \cdot \frac{y-x}{|y-x|}$ on $\O \times \O \setminus \D_{\eta_1}$. We continue to use $C_1 = \frac{1}{2} \big( \frac{ A^2}{\th_0 \th_1} + A \big)$ and $\d'' = \min\{\d', \frac{\th_0^2}{4C_1} \}$ but with the additional assumption that $\d'' << \eta_1$.

\begin{lem}
\label{lem:second-tech}
For $\d < \d''$ sufficiently small and $x, y  \in \O \setminus \O_\d$ with $|x-y| > \eta_1$ there are constants $C_2, C_3$ independent of $\d$ such that:
$$
(X(y) - X(x) ) \cdot \frac{y-x}{|y-x|} \geq \frac{C_2}{\d} - C_3
$$
\end{lem}

\begin{pf}
Using the strict convexity of $\p \O$, choose $\d'' << \eta_1$ so that if $x,y \in \O \setminus \O_{\d''}$ with $|y-x| > \eta_1$ then the line segment  $\g(s)$ joining $x$ to $y$ intersects $\O_{\d''}$.
Suppose $\d < \d''$, $x \in \p \O_{\frac{\d}{2}}$ and $y \in \O$ with $|y-x| > \eta_1$ where we assume that $u_0(x) \leq u_0(y)$. Let $\g(s)$ be the line segment joining $x$ to $y$, parameterized by arc length.  
Divide $\g(s)$ into two disjoint curves:  $\g_1$ lying in $\O \setminus \O_{\d''}$ and $\g_2$ lying in $\O_{\d''}$. Then:
\begin{eqnarray} \nonumber
(X(y) - X(x) ) \cdot \frac{y-x}{|y-x|} & = & \langle X(\g(s)), \g'(s) \rangle|_0^{|y-x|} \\ \nonumber
& = & \int_0^{|y-x|} \frac{d}{ds} ( \langle X(\g(s)), \g'(s) \rangle) ds\\ \nonumber
& = & \int_0^{|y-x|} \n^2(-\log u_0)( \g'(s), \g'(s) ) ds\\ \nonumber
\label{equ:integral-est-1}
& = & \int_{\g_1} \n^2(-\log u_0)( \g'(s), \g'(s) ) ds +  \int_{\g_2} \n^2(-\log u_0)( \g'(s), \g'(s) ) ds\\ 
\end{eqnarray}
From 
$$
|\n^2 \log u_0| \leq \frac{|\n^2 u_0|}{|u_0|} + \frac{|\n u_0|^2}{{u_0}^2} 
$$ 
and (\ref{equ:estimate-A}) we get the estimate:
\begin{equation}
\label{equ:estimate-g1}
\int_{\g_2} \n^2(-\log u_0)( \g'(s), \g'(s) ) ds > -\big( \frac{A}{2 \d''} +  \frac{A^2}{4 {\d''}^2} \big)D
\end{equation}
To estimate the other integral, set $\g'(s) = W$. At each point $\g(s)$, decompose $W$ into a component, $W^{\top}$, tangent to $T_{\g(s)} \O_{u_0(\g(s))}$ and a normal component $W^\perp$ with respect to the inward pointing normal $-\nu_{\g(s)}$. For the curve $\g_1$ lying in $\O \setminus \O_{\d''}$ we have, as in the proof of Lemma \ref{lem:first-tech}:
\begin{eqnarray*}
\n^2 u_0(W,W) & = & \n^2 u_0(W^{\top},W^{\top}) + 2\n^2 u_0(W^{\top},W^\perp)+ \n^2 u_0(W^\perp,W^\perp)\\
& \leq & - \frac{\th_0 \th_1}{2} |W^{\top}|^2 + C_1 |W^\perp|^2\\
\end{eqnarray*}
For $\d < \d''$ let $k$ be the integer such that
$$
2^k \d \leq \d'' < 2^{k+1} \d.
$$
Set $\d_j = 2^j \d$. For $\d_{j-1} \leq u_0 \leq \d_{j}$ with  $j= 0, \dots, k$, we have, as in the proof of Lemma \ref{lem:first-tech}:
\begin{eqnarray}  \nonumber
\n^2 \log u_0(W,W) & = & \frac{\n^2 u_0(W,W)}{u_0}  - \frac{|\n u_0|^2}{u_0^2} |W^\perp|^2 \\ 
\label{equ:inequality-on-path}
& \leq &- \frac{\th_0 \th_1}{2 \d_{j}}|W^{\top}|^2  -  \frac{\th_0^2}{2{\d_{j}}^2}|W^\perp|^2 \\ \nonumber
\end{eqnarray}
We subdivide the curve $\g_1$ via the level sets of $u_0$. Let $s_j$ be the first $s$ satisfying $u_0(\g(s)) =  \d_j$, for $j=-1, 0, \dots, k$. Let $s'$ be the first $s$ satisfying $u_0(\g(s)) =  \d''$. Then $s' > s_k$ and we can write:
\begin{eqnarray*}
\int_{\g_1} \n^2(-\log u_0)( \g'(s), \g'(s) ) ds &\geq& \sum_{j=-1}^{k-1} \int_{s_j}^{s_{j+1}}  \n^2(-\log u_0)( \g'(s), \g'(s) ) ds\\
 &+& \int_{s_k}^{s'}  \n^2(-\log u_0)( \g'(s), \g'(s) ) ds\\ 
 &\geq&\sum_{j=-1}^{k-1} \int_{s_j}^{s_{j+1}}  \n^2(-\log u_0)( \g'(s), \g'(s) ) ds,\\
\end{eqnarray*}
where the last inequality follows from (\ref{equ:integral-positive}). Since $|\n u_0| \leq \frac{A}{2}$ it follows that $s_{j+1} - s_j \geq \frac{\d_{j+1}}{A}$. Thus $s_{j+1}  \geq \frac{\d_{j+1}}{A}$. Similarly $s' \geq \frac{\d''}{A}$. Hence $|y-x| \geq \frac{\d''}{A}$. Since $\O$ is strictly convex, for $\d \leq \d''$, if $x \in \p \O_{\frac{\d}{2}}$ and $y \in  \O$ satisfy $|y-x| \geq \frac{\d''}{A} > 0$ then there exists a constant $\th_2 > 0$, depending only on $ \frac{\d''}{A}$, $\d'$ and the convexity of $\O$ such that:
\begin{equation}
\label{equ:estimate-angle}
\langle -\nu_x, \frac{y-x}{|y-x|} \rangle \geq \th_2.
\end{equation}
This estimate can also be written $|W^\perp| \geq  \th_2$. Therefore using (\ref{equ:inequality-on-path}):
\begin{eqnarray} \nonumber
\sum_{j=-1}^{k-1} \int_{s_j}^{s_{j+1}}  \n^2(-\log u_0)( \g'(s), \g'(s) ) ds &\geq& \sum_{j=-1}^{k-1} (s_{j+1} - s_j)  \frac{\th_0^2}{2{\d_{j}}^2}|W^\perp|^2 \\  \nonumber
\label{equ:estimate-g2}
&\geq& \sum_{j=-1}^{k-1}  \frac{\th_0^2 \th_2^2}{A{\d_{j}}} \\
&\geq& \frac{\th_0^2 \th_2^2}{A \d} \sum_{j=-1}^{k-1} \frac{1}{2^j}  \\  \nonumber
\end{eqnarray}
Set $C_2 = \frac{2 \th_0^2 \th_2^2}{A}$  and $C_3 = \big( \frac{A}{2 \d''} +  \frac{A^2}{4 {\d''}^2} \big)D$ then the estimates (\ref{equ:estimate-g1}) and (\ref{equ:estimate-g2}) imply:
\begin{equation}
(X(y) - X(x) ) \cdot \frac{y-x}{|y-x|} \geq \frac{C_2}{\d} - C_3
\end{equation}

\end{pf}

\bigskip

\begin{thm} 
\label{thm:final-est-on-C}
On $\O \times \O$:
$$
{\cal C}(x,y) = \big( Y(y) - Y(x) \big) \cdot \frac{y-x}{|y-x|} + 2 \psi \bigg( \frac{|y-x|}{2} \bigg) \geq 0
$$
\end{thm}

\begin{pf}
Given $\e > 0$ by Lemma \ref{lem:first-tech} there is an $\eta_1 > 0$ such that ${\cal C}(x,y) \geq -\e$ for $x,y \in \D_{\eta_1}$. By Lemma \ref{lem:second-tech} there is a $\d > 0$ such that for $x \in (\O \setminus \O_\d) \setminus \D_{\eta_1}$ and for any $y \in \O  \setminus \D_{\eta_1}$, ${\cal C}(x,y) \geq 0$. By Theorem \ref{thm:estimate-on-C} this implies that ${\cal C}(x,y) \geq -\e$ for $x, y \in \O$. Since $\e$ is arbitrary this implies the result for $\psi$ defined on $[0, \frac{D'}{2})$. Let $D' \to D$ to conclude the result for $\psi$ satisfying (\ref{equ:psi-neg}) and (\ref{equ:psi-eq1}) on $[0, \frac{D}{2})$.
\end{pf}

Under the assumption that there is a $C^2$ function $\psi$ on $[0, \frac{D}{2})$ that satisfies (\ref{equ:psi-neg}) and (\ref{equ:psi-eq1}), this result implies that $-\psi = \o$ is a modulus of expansion of $Y$.

\bigskip

\subsection{Differential inequalities and a Sturm-Liouville problem}

To apply Theorem \ref{thm:final-est-on-C} we must find a solution to the differential inequalities (\ref{equ:psi-neg}) and (\ref{equ:psi-eq1}).
The inequality (\ref{equ:psi-eq1}) becomes:
\begin{equation}
\label{equ:psi-eq2}
 2\psi''(s) + 4\psi'(s)  \psi(s) + 2\L - 2\t(s) \leq 0
\end{equation}
To proceed we consider two cases:

\medskip

\noindent (I) $\t(s) \geq 0$. This is a ``convexity'' condition on the vector field $V$. In this case (\ref{equ:psi-eq2}) follows from:
\begin{equation}
\label{equ:psi-eq3}
 2\psi''(s) + 4\psi'(s)  \psi(s) + 2\L  \leq 0
\end{equation}

\medskip

\noindent (II) In general, let
$$
\L' = \max \big( \sup_{s \in [0, \frac{D}{2}]} - \t(s) , 0 \big)
$$
and set $\tilde{\L} = \L + \L'$. Then (\ref{equ:psi-eq2}) follows from:
\begin{equation}
\label{equ:psi-eq4}
 2\psi''(s) + 4\psi'(s)  \psi(s) + 2 \tilde{\L}  \leq 0
\end{equation}
In both cases the inequality has the same form. We will set $\L = \tilde{\L}$ and use the inequality:
\begin{equation}
\label{equ:psi-eq5}
 2\psi''(s) + 4\psi'(s)  \psi(s) + 2\L  \leq 0, \;\; \mbox{on} \;\; [0, \tfrac{D}{2}).
\end{equation}
If $g(s)$ is a continuous piecewise differentiable function on $[0, \frac{D}{2}]$ such that $g'(s) \geq \L$ for $s \in [0, \frac{D}{2}]$ then this inequality follows from:
\begin{equation}
\label{equ:psi-eq6}
 \psi'(s) +   \psi(s)^2 + g(s)   = -\nu, \;\; \mbox{on} \;\; [0, \tfrac{D}{2}).
\end{equation}
by differentiation, where $\nu$ is an arbitrary constant.

Make the substitution $\o = -\psi$. Then (\ref{equ:psi-eq6}) becomes:
\begin{equation}
\label{equ:psi-eq7}
\o'(s) - \o(s)^2 - g(s)  = \nu, \;\; \mbox{on} \;\; [0, \tfrac{D}{2}).
\end{equation}
This is a Riccati equation. We will show that this equation can be solved for suitable $\nu$ with $\o'(s) > 0$ for $s \in [0, \frac{D}{2})$ and $\o(0) \leq 0$.

\bigskip

Let $F_{\s}(s)$ be the continuous piecewise differentiable  function:
\begin{equation}
\label{equ:def-F}
F_{\s}(s) = \left\{ \begin{array}{ll}
         \s s & \mbox{if $0 \leq s \leq \tfrac{D}{2}$},\\[.2cm]
         -\s s & \mbox{if $ -\tfrac{D}{2}  \leq s \leq 0$}.\\ \end{array} \right.
\end{equation}

Consider the Sturm-Liouville eigenvalue problem on $[- \frac{D}{2}, \frac{D}{2}]$:
\begin{equation}
\label{equ:SLeigen}
w'' + F_{\s} w = - \mu w
\end{equation}
with $w(- \frac{D}{2}) = w( \frac{D}{2}) = 0$. This is a regular Sturm-Liouville eigenvalue problem in normal form. It has an infinite sequence of real eigenvalues $\mu_0 < \mu_1 < \mu_2 < \dots$ with $\lim_{n \to \infty} \mu_n = \infty$. The eigenfunction ${w}_n(s)$ belonging to the eigenvalue $\mu_n$ has exactly $n$ zeros in the interval $ (-\frac{D}{2}, \frac{D}{2}) $ and is uniquely determined up to a constant factor [BR]. 

There is a scalar $\s_0 > 0$ such that the linear operator on smooth functions on $[- \frac{D}{2}, \frac{D}{2}]$:
$$
{\cal L}(w) = - w'' - F_{\s} w
$$
is  positive definite for $\s < \s_0$ and is not positive definite for $\s > \s_0$. Therefore the first eigenvalue $\mu_0(\s)$ is positive if $\s < \s_0$, zero if $\s = \s_0$ and negative otherwise.  In all cases the corresponding eigenfunction $w_0$ is positive on $(- \frac{D}{2}, \frac{D}{2})$, vanishes at the endpoints and satisfies:
\begin{equation}
\label{equ:convex-w}
w_0''  = -( \mu_0 + F_{\s}(s)) w_0
\end{equation}
Using the variational characterization of the first eigenvalue it follows that:
\begin{enumerate}
\renewcommand{\labelenumi}{(\roman{enumi})}
\item $-\mu_0 = -\mu_0(\s)$ is an increasing function of $\s$,\\
\item $-\mu_0 < \s \tfrac{D}{2}$.
\end{enumerate}
Since $F_\s$ is an even function, so is $w_0$ and therefore $w_0'(0)=0$. Thus $w_0$ satisfies the boundary value problem on $[0, \frac{D}{2}]$:
\begin{eqnarray}
\label{BVP-w}
&&w''(s) + F_\s(s) w(s) =  - \mu_0 w(s) \\ \nonumber
&& w'(0) = 0, w(\tfrac{D}{2}) =0.
\end{eqnarray}

\begin{prop}
\label{prop:comparison}
As $\s \to \infty$, $\frac{-\mu_0}{\s} \to  \frac{D}{2}$.
\end{prop}

\begin{pf}
Set $y(s) = w_0(s \s^{-\frac{1}{3}})$. Then,
\begin{eqnarray*}
y''(s) &=& \s^{-\frac{2}{3}} w_0''(s \s^{-\frac{1}{3}}) \\
&=& \s^{-\frac{2}{3}} \big( -\s s \s^{-\frac{1}{3}} w_0(s \s^{-\frac{1}{3}}) - \mu_0 w_0(s \s^{-\frac{1}{3}}) \big)\\
&=& - s y(s) - \mu_0 \s^{-\frac{2}{3}} y(s).
\end{eqnarray*}
Hence $y$ is the first eigenfunction with eigenvalue $ -\mu_0  \s^{-\frac{2}{3}}$ of the boundary value problem on $[0, \tfrac{D}{2} \s^{\frac{1}{3}}]$:
\begin{eqnarray}
\label{BVP-y}
&&y''(s) + s y(s) = -\mu_0  \s^{-\frac{2}{3}} y(s) \\ \nonumber
&& y'(0) = 0, \;\; y(\tfrac{D}{2} \s^{\frac{1}{3}}) =0.
\end{eqnarray}
Scaling we can suppose that $y$ satisfies: $\int_0^{\frac{D}{2} \s^{\frac{1}{3}}} y(s)^2 ds = 1$. Set $\b = \mu_0  \s^{-\frac{2}{3}}$. Then $\b$ is characterized as the infimum
$$
\b = \inf \big( \int_0^{\frac{D}{2} \s^{\frac{1}{3}}} \big( (y'(s))^2 - s (y(s))^2 \big) ds \big)
$$ 
over functions $y$ satisfying $\int_0^{ \frac{D}{2}\s^{\frac{1}{3}}} y(s)^2 ds = 1$. Clearly $\b > - \s^{\frac{1}{3}} \frac{D}{2}$. On the other hand define the test function:
\begin{equation}
z(s) = \left\{ \begin{array}{ll}
         \sqrt{2} \sin \pi(s-  \frac{D}{2} \s^{\frac{1}{3}})  & \mbox{if $ \frac{D}{2} \s^{\frac{1}{3}} - 1 \leq s \leq  \frac{D}{2} \s^{\frac{1}{3}} $},\\[.2cm]
         0 & \mbox{if $ 0  \leq s \leq  \frac{D}{2} \s^{\frac{1}{3}}  - 1$}.\\ \end{array} \right.
\end{equation}        
Then,
$$
\int_0^{\frac{D}{2}\s^{\frac{1}{3}}} z(s)^2 ds = 1
$$
and 
$$
\int_0^{\frac{D}{2} \s^{\frac{1}{3}}} \big( (z'(s))^2 - s (z(s))^2 \big) ds =- \tfrac{D}{2} \s^{\frac{1}{3}} + \pi^2 + \tfrac{1}{2}.
$$
Hence,
$$
 -  \tfrac{D}{2} \s^{\frac{1}{3}} < \b  \leq - \tfrac{D}{2} \s^{\frac{1}{3}} + \pi^2 + \tfrac{1}{2}.
$$
Since $\b = \mu_0  \s^{-\frac{2}{3}}$ this implies,
$$
\tfrac{D}{2} - (\pi^2 + \tfrac{1}{2}) \s^{-\frac{1}{3}} \leq   \frac{- \mu_0}{\s} < \tfrac{D}{2}.
$$
The result follows.
\end{pf}

Introduce the function
$$
v_\s(s) = \frac{w_0'(s)}{w_0(s)} \;\;\;\; 0 \leq s \leq \frac{D}{2}.
$$
On $[0, \frac{D}{2}]$, $v_\s(s)$ satisfies the initial value problem:
\begin{eqnarray}
\label{IVP-v}
&&v'(s) + v(s)^2 = - \s s - \mu_0 \\ \nonumber
&& v(0) = 0
\end{eqnarray}

\begin{prop}
For each $\s > \s_0$, there is  a unique point $s_0= s_0(\s) \in (0, \frac{D}{2})$ such that $v_\s'(s_0) = 0$, $v_\s'(s) > 0$ on $(0, s_0)$ and $v_\s'(s) < 0$ on $(s_0, \frac{D}{2})$. Moreover, as $\s \to \infty$, $s_0(\s) \to 0$. 
\end{prop}

\begin{pf}
Differentiating (\ref{IVP-v}) we have:
$$
v''(s) = -\s - 2 v(s) v'(s)
$$
Therefore every critical point of $v$ in $(0, \frac{D}{2})$ is a local maximum. Note that provided $\mu_0 < 0$, $v'(0) > 0$ so that $v$ is initially increasing and positive. On the other hand, since $0 < - \mu_0 < \s$, by (\ref{equ:convex-w})  there is an  $s_1 \in (0, \frac{D}{2})$ such that $w''(s) < 0$ for $s \in (s_1, \frac{D}{2})$. Thus on  $(s_1, \frac{D}{2})$, $v'(s) = \frac{w''(s)}{w(s)} - (\frac{w'(s)}{w(s)})^2 < 0$. It follows that $v$ has at least one local maximum point in $(0,  \frac{D}{2})$. Thus there is a unique critical point $s_0 \in (0, \frac{D}{2})$ and it is a local maximum. Therefore $v'(s) > 0$ on $[0, s_0)$ and $v'(s) < 0$ on $(s_0, \frac{D}{2})$. From (\ref{IVP-v}) it follows at that the maximum point $s_0$:
$$
v(s_0)^2 = -\s s_0 - \mu_0
$$
Set $\l = -\mu_0$. Thus
$$
\max_{[0, \frac{D}{2}]} v(s) = v(s_0) = \sqrt{\l - \s s_0} < \sqrt{\l}
$$

Let $0 < a < \frac{D}{2}$. We will show that for $\s$ sufficiently large (depending on $\frac{1}{a^2}$) the point $s_0 \in (0,a)$. Suppose not. Then on $[0, a)$, $v'(s) > 0$ and $v$ is strictly increasing. We deduce a contradiction by showing that then $v(a) > \max_{[0, \frac{D}{2}]} v(s)$. Let $0 = t_0 < t_1 < \dots < t_n =a$ be a subdivision of $[0, a]$ with $t_{i+1} - t_i = \D t$. Denote $v(t_i) = v_i$. We construct an iterative scheme that successively estimates $v_i$ above and below. We denote the lower estimate $\bar{v}_i$, the upper estimate ${\bar{\bar{v}}}_i$ so that $\bar{v}_i \leq v_i \leq {\bar{\bar{v}}}_i$.
At $t_0 = 0$ set $\bar{v}_0 = v_0 = {\bar{\bar{v}}}_0 = 0$. On the interval $[t_0, t_1]$ we have:
$$
-\s t_1 + \l - v_1^2 \leq v' \leq \l.
$$
Hence,
$$
v_1 \leq \l \D t,
$$
so that,
$$
-\s t_1 + \l - (\l \D t)^2 \leq v' \leq \l
$$
Therefore,
\begin{eqnarray*}
{\bar{\bar{v}}}_1 &=& \l  \D t,\\
\bar{v}_1 &=& (-\s t_1 + \l - (\l \D t)^2) \D t.
\end{eqnarray*}
Suppose that $\bar{v}_i$ and ${\bar{\bar{v}}}_i$ are defined. Then on $[t_i, t_{i+1}]$ we have:
$$
-\s t_{i+1} + \l - {\bar{\bar{v}}}_{i+1}^2 \leq v' \leq -\s t_{i} + \l - \bar{v}_i^2.
$$
Hence we define:
$$
{\bar{\bar{v}}}_{i+1} = (-\s t_{i} + \l - \bar{v}_{i}^2) \D t + {\bar{\bar{v}}}_i, 
$$
so that, $v_{i+1} \leq {\bar{\bar{v}}}_{i}$. Then on $[t_i, t_{i+1}]$:
$$
-\s t_{i+1} + \l -\big( (-\s t_{i} + \l - \bar{v}_i^2) \D t + {\bar{\bar{v}}}_i \big)^2 \leq v'.
$$
Define,
$$
\bar{v}_{i+1} = \bigg(-\s t_{i+1} + \l -\big( (-\s t_{i} + \l - \bar{v}_i^2) \D t + {\bar{\bar{v}}}_i \big)^2 \bigg) \D t + \bar{v}_i.
$$
Thus since,
$$
{\bar{\bar{v}}}_{i} = (-\s t_{i-1} + \l - \bar{v}_{i-1}^2) \D t + {\bar{\bar{v}}}_{i-1}
$$
we derive,
\begin{eqnarray*}
\bar{v}_{i+1} &\equiv& (-\s t_{i+1} + \l - {\bar{\bar{v}}}_i^2) \D t +  \bar{v}_i \;\; \mbox{mod} (\D t)^2, \\
&\equiv& (-\s t_{i+1} + \l - {\bar{\bar{v}}}_{i-1}^2) \D t +  \bar{v}_i \;\; \mbox{mod}  (\D t)^2,\\
&\equiv& (-\s t_{i+1} + \l - {\bar{\bar{v}}}_{i-2}^2) \D t +  \bar{v}_i \;\; \mbox{mod}  (\D t)^2,\\
& \dots & \\
&\equiv& (-\s t_{i+1} + \l ) \D t +  \bar{v}_i \;\; \mbox{mod}  (\D t)^2.
\end{eqnarray*}
Hence,
$$
v(a) \geq \bar{v}_n \equiv (-\s t_n + \l) \D t +  (-\s t_{n-1} + \l) \D t + \dots + \l \D t \;\; \mbox{mod} (\D t)^2.
$$
Letting $n \to \infty$ we get,
$$
v(a) \geq  -\s \int_0^a t dt + \l a = -\s \frac{a^2}{2} + \l a.
$$
Recall that,
$$
v(a) \leq \max_{[0, \frac{D}{2}]} v \leq \sqrt{\l}.
$$
Hence,
$$
-\s \frac{a^2}{2} + \l a <  \sqrt{\l}.
$$
Using Proposition \ref{prop:comparison} there is a scalar $\s_1$ such that if $\s \geq \s_1$ then $\s \frac{D}{3} < \l < \s \frac{D}{2}$. Hence:
$$
a(-\s \frac{a}{2} +\s \frac{D}{3}) < -\s \frac{a^2}{2} + \l a <  \sqrt{\l} < \sqrt{\s \frac{D}{2}}
$$
Since $a < \frac{D}{2}$ we get,
$$
 \s  <   \frac{72}{D} \frac{1}{a^2}.
$$
Therefore, if,
$$
 \s  \geq \max(\s_1,  \frac{72}{D} \frac{1}{a^2})
$$
we have a contradiction, proving that, in this case, the maximum of $v$ occurs in $(0, a)$.
\end{pf}

In particular, there is a $\s_2 > \s_0$ such that for all $\s > \s_2$ we have $s_0(\s) \in (0, \frac{D}{4})$. 


Set:
$$
\eta_\s = (\tfrac{D}{2})^{-1} ( \tfrac{D}{2} - s_0(\s))
$$
and note that for $\s > \s_2$, $\eta_\s > \tfrac{1}{2}$, independent of the choice of $\s$. Let $\eta$ satisfy:
\begin{equation}
\label{equ:scaling-parameter}
\eta_\s \geq \eta > \tfrac{1}{2}
\end{equation}
We next explain how to use this Sturm-Liouville problem to solve (\ref{equ:psi-eq7}).

For $s \in [0, \frac{D}{2})$ and $\eta$ satisfying (\ref{equ:scaling-parameter}) set:
$$
\o(s) = -\frac{ \eta w_0'(\eta s + s_0)}{w_0(\eta s +s_0)}.
$$ 
Set
$$
\tilde{F}_\s(s) = \eta^2 F_\s( \eta s + s_0)
$$
Then
\begin{eqnarray}
\nonumber
\o'(s ) - \o^2(s ) - \tilde{F}_\s(s)  &=&  -\eta^2 \bigg( \frac{w_0''(\eta s + s_0)}{w_0(\eta s + s_0)} + F_\s(\eta s + s_0) \bigg)\\ \nonumber
& = & -\eta^2 \bigg(\frac{1}{w_0(\eta s + s_0)} \big( w_0''(\eta s + s_0) + F_\s(\eta s + s_0) w_0(\eta s + s_0) \big)\bigg) \\ \label{equ:solving-riccatti}
&= & \eta^2 \mu_0
\end{eqnarray}
Hence if we chose $\nu = \eta^2 \mu_0$ then $\o(s)$ satisfies  (\ref{equ:psi-eq7}) with $g(s) = \tilde{F}_\s(s)$. Since,
$$
 \frac{d}{ds}\tilde{F}_\s(s) = \eta^3 F_\s'( \eta s + s_0) = \eta^3 \s
 $$
Choose
\begin{equation}
\label{equ:choice-of-sigma}
\s = \max( \s_2, 8 \L),
\end{equation} 
 then $g'(s) \geq \L$ and (\ref{equ:psi-eq7}) follows. Note that for $s \in [0, \frac{D}{2})$:
 $$
 \o(s) = -\eta v_\s(\eta s + s_0)
$$
Therefore  for $s \in [0, \frac{D}{2})$:
$$ 
  \o'(s) = -\eta^2 v'_\s(\eta s + s_0) > 0
 $$
 Clearly, $\o(0) =  -\eta v_\s( s_0) < 0$. It follows that $\psi(s) = -\o(s)$ on $[0, \frac{D}{2})$ satisfies
the requirements of Theorem \ref{thm:estimate-on-C}.

Notice that if we choose the scaling factor $\eta$ satisfying $\eta_\s > \eta > \frac{1}{2}$ then $\o(s)$ is defined on some interval $[0, \frac{D'}{2})$ with $D' > D$. 
In particular $\o(s)$ is defined on $[0, \frac{D}{2}]$. This choice is used in the next subsection.

\bigskip

\subsection{Associated Sturm-Liouville problem}

In this subsection we show that the Sturm-Liouville  problem: 
\begin{equation}
\label{equ:ST-for-w2}
{\cal L} [{w}] = {w}''+ F_\s(s)  {w} = -\mu {w}
\end{equation}
with $w( \pm \frac{D}{2}) = 0$, where $F_\s$ is defined above, is an associated Sturm-Liouville problem and thus yields a spectral gap result.

The second eigenfunction $w_1$ has a unique zero and because $F_\s$ is even this zero occurs at the origin.  After scaling we can assume that $w_0(s) > 0$ for $s \in (-\frac{D}{2}, \frac{D}{2})$  and $w_1(s) > 0$ for $s \in (0, \frac{D}{2})$. 

Given any two smooth functions $u,v$ on $[-\frac{D}{2} , \frac{D}{2}]$ the Lagrange identity is easily derived:
$$
u {\cal L}[v] - v {\cal L}[u]  = \frac{d}{ds} \big(u(s) v'(s) - v(s) u'(s) \big).
$$
Set $u = w_1$ and $v = w_0$, then the Lagrange identity yields:
\begin{equation}
\label{equ:lagrange-applic}
w_1 {\cal L}[ w_0] -  w_0 {\cal L}[w_1]  = \frac{d}{ds} \big(w_1  w'_0 -  w_0 w'_1 \big).
\end{equation}

By our scaling assumptions $\int_{0}^{\frac{D}{2}} w_1 w_0  > 0$. Since $\int_{-\frac{D}{2}}^{\frac{D}{2}} w_1 w_0  = 0$ and $w_0 w_1 < 0$ on $(-\frac{D}{2}, 0)$, it follows that  $\int_{s}^{\frac{D}{2}} w_1 w_0 > 0$ for every $s > -\frac{D}{2}$. Integrating (\ref{equ:lagrange-applic}) from $s$ to $\frac{D}{2}$ we have:
$$
0 < - ( \mu_0 - \mu_1) \int_s^{\frac{D}{2}} w_0 w_1  = -\big(w_1(s)  w'_0(s) -  w_0(s) w'_1(s) \big),
$$
where we have used that $w_0(\frac{D}{2}) = w_1(\frac{D}{2}) = 0$. Hence, for any $s \in (-\frac{D}{2}, \frac{D}{2})$:
\begin{equation}
\label{equ:inequality-on-eigen}
\big(w_0(s)  w'_1(s) -  w_1(s) w'_0(s) \big) > 0.
\end{equation}
In particular, this implies that for $s \in (-\frac{D}{2}, \frac{D}{2})$ we have:
\begin{equation}
\label{equ:inequality-on-deriv}
\frac{d}{ds} \big( \frac{w_1}{w_0} \big) (s) = \frac{1}{w_0^2}(s)\big(w_0(s)  w'_1(s) -  w_1(s) w'_0(s) \big) > 0 
\end{equation}

\medskip

Set $\r(s, t) = \frac{e^{-\mu_1 t} w_1(s)}{ e^{-\mu_0 t} w_0(s)}$. Denote $\frac{\p}{\p s} \r(s,t) = \r'(s,t)$. Then 

\begin{lem}
For $s \in (-\frac{D}{2}, \frac{D}{2})$, 
$$ \frac{\p}{\p t} \r(s,t)  = -(\mu_1-\mu_0) \r(s,t) = \r''(s, t) + 2 \frac{w_0'(s)}{w_0(s)} \r'(s, t).$$
\end{lem}

\begin{pf} This is a direct computation.
\end{pf}

To exploit the lemma, choose the scaling factor $\eta$  to satisfy $\eta_\s > \eta > \frac{1}{2}$ and the point $s_0(\s)$ as defined above. Set:
\begin{equation}
\label{eq:varphi}
\varphi(s,t) = C e^{-(\mu_1 - \mu_0) \eta^2 t} \frac{w_1 (\eta s + s_0)}{w_0(\eta s + s_0)},
\end{equation}
where $C$ is a constant to be determined. Note that since $\eta$ satisfies $\eta_\s > \eta > \frac{1}{2}$, $\varphi(s,t)$  is defined on some interval $[0, \frac{D'}{2})$ with $D' > D$. Hence $\varphi(s,t)$ is defined on  $[0,  \frac{D}{2}]$.  Recall that: $\o(s) = -\frac{ \eta w_0'(\eta s + s_0)}{w_0(\eta s + s_0)}$. By the same computation as in the lemma  we get:
$$
\frac{\p \varphi}{\p t} = \varphi'' - 2 \o \varphi' \;\; \mbox{ on} \;\;   [0, \frac{D}{2}] \times \R_+,
$$
Set  
$$
\varphi_0(s) = \varphi(s, 0) = C \frac{w_1(\eta s + s_0)}{w_0(\eta s + s_0)}.
$$ 
Since $w_0(s_0) > 0$ and $w_1(s_0) > 0$ it follows that both $\varphi_0(0) >  0$ and $\varphi(0, t) > 0$ for all $t > 0$. It is also true
that $\frac{\p}{\p s} \varphi(s,t) = C e^{-(\mu_1 - \mu_0) \eta^2 t} \frac{d}{ds} (\frac{w_1}{w_0})(\eta s + s_0) > 0$ for all $t \geq 0$ and $0 \leq s \leq \frac{D}{2}$.

We prove:

\begin{thm}
Let $\O$ be a strictly convex bounded domain with smooth boundary in $\R^n$. Suppose the diameter of $\O$ is $D$. Then the gap between the first   eigenvalue, $\l_0$, and the real part of any other eigenvalue, $\l$, of the linear elliptic operator $L$ given by (\ref{equ:def-of-L}) satisfies:
$$
\Real(\l) - \l_0 > \frac{1}{4} (\mu_1 - \mu_0) = \a > 0,
$$
where $\a$ is a constant depending on ${\L}$ and hence on $b^i, c, \k$.
\end{thm}

\begin{pf}
We derive the theorem from Theorem \ref{thm:basicAC}. Suppose that the eigenfunction corresponding to $\l$ is denoted $u$. Let $z(x,t) = e^{-(\l - \l_0)t} \frac{u(x)}{u_0(x)}$. 
We wish to conclude that for suitable constant $C$, $\varphi(s,t)$ is a modulus of continuity of $z(x,t)$.

The drift velocity $Y( \cdot, t)$ has modulus of expansion 
$$
\o(\cdot, t) = -\psi(\cdot, t) = -(\log w_0(\eta  s + s_0))' =-( \frac{\eta  w_0'}{w_0})(\eta  s + s_0),
$$
by Theorem \ref{thm:final-est-on-C}. As shown above the function $\varphi(s, t)$ satisfies the equation:
$$
\frac{\p \varphi}{\p t} \geq \varphi'' - 2 \o \varphi' \;\; \mbox{ on} \;\;   [0, \frac{D}{2}] \times \R_+,
$$

Using (\ref{equ:inequality-on-deriv}) $\frac{d}{ds} \big(  \frac{w_1}{w_0}(\eta  s + s_0) \big) > 0$ on $[0, \frac{D}{2}]$. Therefore there exists a constant C such that $\frac{d}{d s} \big(C \frac{w_1}{w_0}(\eta s + s_0) \big)$ is a modulus of continuity of $z(x,0)$. Set
$$\varphi(s,t)= C e^{-(\mu_1 - \mu_0) \eta^2 t} ( \frac{w_1}{w_0})(\eta  s + s_0) ,$$
with this constant $C$. The hypotheses of Theorem \ref{thm:basicAC} are satisfied and we conclude that
 $2 \varphi(s,t)$ is a modulus of continuity of $z(x,t)$. Thus,
 $$
  |e^{-(\l -\l_0)t}| \bigg| \frac{u(y)}{u_0(y)} -   \frac{u(x)}{u_0(x)} \bigg|  \leq  C' e^{-(\mu_1 - \mu_0) \eta^2 t}  \big(\frac{w_1}{w_0} \big)(\eta \tfrac{|y-x|}{2} + s_0)  ,
 $$
for some constant $C'$ and any $t \geq 0$. Therefore $\Real(\l) - \l_0 \geq \eta^2 (\mu_1 - \mu_0) > \frac{1}{4} (\mu_1 - \mu_0)$. 
 \end{pf}

\bigskip

\section{\bf Special case: The $\phi$-Laplacian}

Let $\O$ be a bounded domain in $\R^n$ with smooth strictly convex boundary and let $\phi$ and $c$ be  smooth functions of $\bar{\O}$. Consider the eigenvalue problem with Dirichlet boundary conditions on $\O$:
\begin{equation}
\label{equ:eigenvalue-pde}
\D u -  \n \phi \cdot \n u - c u = -\l u.
\end{equation}
Introduce the $\phi$-Laplacian $\D_\phi = \D - \n \phi \cdot \n$, also called the Bakry-Emery Laplacian,  to rewrite  (\ref{equ:eigenvalue-pde}) as:
\begin{equation}
\label{equ:eigenvalue-pde-2}
\D_\phi u - c u = -\l u.
\end{equation}
with $u = 0$ on $\p \O$. The operator $\D_\phi  - c$ is not $L^2$-symmetric with respect to the euclidean volume form $dv$ however introducing the weighted volume form $e^{-\phi} dv$ it is easy to show that this operator is symmetric with respect to the $L^2$ inner product with Dirichlet boundary conditions. Hence from $L^2$ elliptic theory [E] the eigenvalues of (\ref{equ:eigenvalue-pde}) are real.

Comparing with the eigenvalue problem $L(u) = -\l u$ with $L$ as in (\ref{equ:def-of-L}) we have:
$$
b^i = \n_{x_i} \phi.
$$
Hence,
$$
U^{ij} =  \n_{x_j} b^i - \n_{x_i} b^j= 0 .
$$
It follows that Proposition \ref{prop:bound-on-beta} is not needed and in Theorem \ref{thm:main-intro} the parameter $\s$ depends only on $|| \phi ||_{C^4(\bar{\O})}$ and $|| c ||_{C^2(\bar{\O})}$ and does not depend on $K$ or on the geometry of $u_0$. Moreover,
$$
 V^j =   \n_{x_j} c + \tfrac{1}{4}  \n_{x_j} (|B|^2) - \tfrac{1}{2} \D b^j= \n_{x_j} \big(c + \tfrac{1}{4}  |\n \phi|^2 - \tfrac{1}{2} \D \phi \big)
 $$

\begin{defn}
A function $c$ is called $\phi$-convex if the function $c - \frac{1}{2} \D \phi + \frac{1}{4} |\n \phi|^2 = c  -  \frac{1}{2} \D_{\frac{1}{2} \phi} \phi$ is convex in the usual sense. In particular, if
$$
( V(y) - V(x)) \cdot \frac{y-x}{|y-x|} \geq 0
$$
In the notation of the previous section this is equivalent to $\t(s) = 0$.
\end{defn}

\bigskip

\begin{thm}
If $\phi$ is any $C^4$ function and  $g$ is $\phi$-convex  then an associated Sturm-Liouville problem to the eigenvalue problem (\ref{equ:eigenvalue-pde}) is:
$$
w'' + \mu w = 0
$$
on $[-\frac{D}{2}, \frac{D}{2}]$ with $w(-\frac{D}{2}) = w(\frac{D}{2}) = 0$.
Hence the spectral gap for the operator $\D_\phi  - c$ on a convex domain $\O$ satisfies:
$$
 { \l}_1 -  { \l}_0 \geq \frac{3 \pi^2}{D^2}.
$$
\end{thm}

In the case that $\phi$ is a constant this is the result of [AC]. 

\bigskip

\section{\bf General remarks}

It is well known that, for example, the Schr\"odingier operator with a  double well potential on $\R^n$ does not satisfy a uniform non zero gap between the first and second 
eigenvalues. Harrell [H] gives a family of such operators with a separation $R$ between the pairs of wells. As $R \to \infty$ the eigenvalue gap goes to zero. Of course,
theses examples do not apply to a bounded domain. However, we note that the  Sturm-Liouville problem used above (\ref{equ:def-F})   and (\ref{equ:SLeigen})
has the property that as $\s \to \infty$ the eigenvalue gap of the the Sturm-Liouville problem goes to zero. This show that the method used here does not yield useful results without suitable bounds.  
It suggests, though does not prove, that as the constant $\L \to \infty$
the eigenvalue gap of the Dirichlet problem for the operator (\ref{equ:def-of-L}) also goes to zero. 

The Sturm-Liouville problem (\ref{equ:def-F})  and (\ref{equ:SLeigen}), while natural for the problem, is somewhat arbitrary. It is not difficult to find other potential functions whose eigenfunctions yield solutions to the differential inequalities (\ref{equ:psi-neg}) and (\ref{equ:psi-eq1}). The choice of $F_\s$ in (\ref{equ:def-F}) was made because, since $F_\s$ is even, the unique zero of the second eigenfunction $w_2$ is at the origin. Control of the location of this zero is necessary to complete the proof. It is likely that other choices of potential functions can also determine associated Sturm-Liouville problems. Perhaps some of these problems give better spectral gap results.

\end{document}